\documentclass{article}

\usepackage{amsfonts,amssymb,amsthm,epic,eepic}
\usepackage[totalwidth=17cm,totalheight=24cm]{geometry}
\usepackage[T1]{fontenc}

\newcommand{\dr}{\partial}
\newcommand{\C}{{\mathbb C}}

\newcommand{\R}{{\mathbb R}}

\newcommand{\gt}{\tilde{g}}
\newcommand{\Ct}{\tilde{C}}

\newcommand{\Mt}{\tilde{M}}
\newcommand{\Nt}{\tilde{N}}

\newtheorem{prop}{Proposition}[section]
\newtheorem{lemma}[prop]{Lemma}
\newtheorem{sublemma}[prop]{Sub-lemma}

\newtheorem{thm}[prop]{Theorem}
\newtheorem{cor}[prop]{Corollary}
\newtheorem{remark}[prop]{Remark}

\newtheorem{df}[prop]{Definition}
\newtheorem{pty}[prop]{Property}
\newtheorem{question}[prop]{Question}

\newcommand{\pg}{\paragraph}

\newenvironment{thn}[1]{\vskip 0.2cm \noindent{\bf Theorem #1.} \it}{\rm
\vspace{0.2cm}} 
\newenvironment{crn}[1]{\vskip 0.2cm \noindent{\bf Corollary #1.} \it}{\rm
\vspace{0.2cm}} 
\newenvironment{lmn}[1]{\vskip 0.2cm \noindent{\bf Lemma #1.} \it}{\rm
\vspace{0.2cm}} 
\newenvironment{qn}[1]{\vskip 0.2cm \noindent{\bf Question #1.} \it}{\rm
\vspace{0.2cm}} 
\newenvironment{sketch}{\vskip 0.2cm \noindent{\bf Brief sketch of the
    proof~. ~~~}}{$\qed$ \vspace{0.2cm}} 

\newcommand{\btm}{\begin{thm}}
\newcommand{\etm}{\end{thm}}
\newcommand{\bpt}{\begin{pty}}
\newcommand{\ept}{\end{pty}}
\newcommand{\blm}{\begin{lemma}}
\newcommand{\elm}{\end{lemma}}
\newcommand{\bsl}{\begin{sublemma}}
\newcommand{\esl}{\end{sublemma}}
\newcommand{\bcr}{\begin{cor}}
\newcommand{\ecr}{\end{cor}}
\newcommand{\bdf}{\begin{df}}
\newcommand{\edf}{\end{df}}
\newcommand{\bprop}{\begin{prop}}
\newcommand{\eprop}{\end{prop}}
\newcommand{\bas}{\begin{asser}}
\newcommand{\eas}{\end{asser}}
\newcommand{\beq}{\begin{equation}}
\newcommand{\eeq}{\end{equation}}
\newcommand{\bpv}{\begin{proof}}
\newcommand{\epv}{\end{proof}}
\newcommand{\bpvs}{\begin{sketch}}
\newcommand{\epvs}{\end{sketch}}
\newcommand{\bit}{\begin{itemize}}
\newcommand{\eit}{\end{itemize}}
\newcommand{\bpn}{\begin{pfn}}
\newcommand{\epn}{\end{pfn}}
\newcommand{\btn}{\begin{thn}}
\newcommand{\etn}{\end{thn}}
\newcommand{\bcn}{\begin{crn}}
\newcommand{\ecn}{\end{crn}}
\newcommand{\bqn}{\begin{qn}}
\newcommand{\eqn}{\end{qn}}
\newcommand{\bln}{\begin{lmn}}
\newcommand{\eln}{\end{lmn}}
\newcommand{\brk}{\begin{remark}}
\newcommand{\erk}{\end{remark}}
\newcommand{\bq}{\begin{question}}
\newcommand{\eq}{\end{question}}

\newenvironment{pfn}[1]{\vskip 0.2cm \noindent{\it Proof #1.}}{$\square$
\vspace{0.2cm}}

\newcommand{\cC}{\mathcal{C}}

\begin{document}

\title{Hyperideal circle patterns}

\author{Jean-Marc Schlenker\thanks{
Laboratoire Emile Picard, UMR CNRS 5580,
Universit{\'e} Paul Sabatier,
118 route de Narbonne,
31062 Toulouse Cedex 4,
France.
\texttt{schlenker@picard.ups-tlse.fr; http://picard.ups-tlse.fr/\~{
}schlenker}. }}

\date{July 2004 (v2)}

\maketitle

\begin{abstract}
A ``hyperideal circle pattern'' in $S^2$ is a finite family of oriented
circles, similar to an ``usual'' circle pattern but such that the closed
disks bounded by the circles do not cover the whole sphere. Hyperideal circle
patterns are directly 
related to hyperideal hyperbolic polyhedra, and also to circle packings. 

To
each hyperideal circle pattern, one can associate an incidence graph and a
set of intersection angles. 
We characterize the possible incidence graphs and intersection angles of
hyperideal circle patterns in the sphere, the torus, and in higher genus
surfaces. It is a consequence of a more general result, describing the
hyperideal circle patterns in the boundaries of geometrically finite
hyperbolic 3-manifolds (for the corresponding $\C P^1$-structures). This more
general statement is obtained as a consequence of a theorem of Otal
\cite{otal,bonahon-otal} on the
pleating laminations of the convex cores of geometrically finite hyperbolic
manifolds. 

\bigskip

\begin{center} {\bf Résumé} \end{center}

%\bigskip

Un ``motif de cercles hyperidéal'' sur $S^2$ est une famille finie de cercles
orientés, similaire à un motif de cercles ``usuel'' mais tel que la réunion
des disques fermés bordés par les cercles ne couvre pas la sphère. Les motifs
de cercles hyperidéaux sont directement liés aux polyèdres hyperboliques
hyperidéaux, et aussi aux empilements de cercles. 

A chaque motif de cercles hyperidéal, on associe un graphe d'incidence et un
ensemble d'angles d'intersection. On caractérise les graphes d'incidence et
les angles d'intersection possibles dans la sphère, le tore, et sur les
surfaces de genre supérieur. C'est une conséquence d'un résultat plus général
décrivant les motifs de cercles hyperidéaux dans les bords de variétés
hyperboliques géométriquement finies (pour les $\C P^1$-structures
correspondantes). Ce résultat plus général est obtenu comme conséquence d'un
théorème d'Otal \cite{otal,bonahon-otal} sur les
laminations de plissage des coeurs convexes de variétés hyperboliques
géométriquement finies.
\end{abstract}

\section{Introduction}

\pg{Circle patterns.}
A {\it circle packing} on the sphere $S^2$ is a finite family of oriented
circles with disjoint interiors. The {\it incidence graph} of a circle packing
is a graph $\Gamma$, embedded in $S^2$, which has one vertex for each circle,
and an edge between two vertices when the corresponding circles are tangent. A
classical theorem of Koebe \cite{koebe} states that, given any graph $\Gamma$
on $S^2$ which is the 1-skeleton of a polytopal
triangulation, there is a unique circle
packing with incidence graph $\Gamma$ --- the uniqueness is up to M\"obius
transformations. This theorem was extended to higher genus surfaces by
Thurston \cite{thurston-notes}. It also holds when $\Gamma$ is the 1-skeleton
of a cellular decomposition; one should then add the condition that, for each
connected component of the complement of the disks, there is another circle
orthogonal to all the adjacent circles of the pattern.

One can also consider {\it circle patterns} on $S^2$ (see
e.g. \cite{bobenko-springborn}).
For reasons that should become apparent below, we add the
adjective {\it ideal} to describe what is perhaps the most commonly considered
type of circle patterns.

\begin{df} \label{df:ideal}
A {\bf circle pattern} on $S^2$ is a finite family of oriented circles
$C_1, \cdots, C_N$. Given a circle pattern, an {\bf interstice} is a connected
component of the complement of the union of the open disks bounded by the
circles. If $C_1, \cdots, C_N$ is a circle pattern, it is {\bf ideal} if:
\begin{itemize}
\item Each interstice is a point. 
\item If $D$ is an open disk in $S^2$, containing no interstice, but such that
  its closure contains at least 3 of the interstices, then $D$ is the open
  disk bounded by one of the $C_i, 1\leq i\leq N$.
\end{itemize}
\end{df}

It follows from this definition that any interstice is in at least $3$ of the
circles. However there might also be points which are in $3$ circles or more,
but which are not interstices. 

Given an ideal circle pattern, its {\it incidence graph} is the 1-skeleton of
the cellular decomposition of the sphere which has:
\begin{itemize}
\item One vertex for each circle.
\item One face for each interstice.
\item One edge between two vertices, when the corresponding circles intersect
  at two interstices. 
\end{itemize}

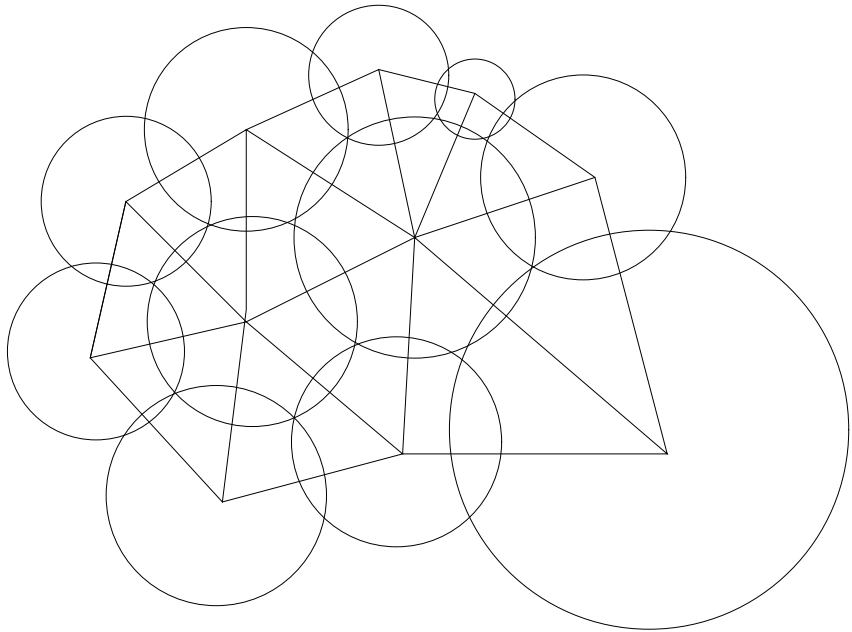
\begin{figure}[htbp]
\begin{center}
\setlength{\unitlength}{0.00083333in}
\begingroup\makeatletter\ifx\SetFigFont\undefined%
\gdef\SetFigFont#1#2#3#4#5{%
  \reset@font\fontsize{#1}{#2pt}%
  \fontfamily{#3}\fontseries{#4}\fontshape{#5}%
  \selectfont}%
\fi\endgroup%
{\renewcommand{\dashlinestretch}{30}
\begin{picture}(5266,3928)(0,-10)
\put(749,2680){\ellipse{1060}{1060}}
\put(1499,3130){\ellipse{1272}{1272}}
\put(1536,1930){\ellipse{1310}{1310}}
\put(2549,2455){\ellipse{1508}{1508}}
\put(2324,3468){\ellipse{876}{876}}
\put(2436,1180){\ellipse{1310}{1310}}
\put(1311,843){\ellipse{1378}{1378}}
\put(561,1743){\ellipse{1106}{1106}}
\put(4011,1255){\ellipse{2494}{2494}}
\put(2924,3318){\ellipse{502}{502}}
\put(3599,2830){\ellipse{1282}{1282}}
\path(1499,3130)(749,2680)(524,1705)
	(1349,805)(1499,2005)(1499,3130)
\path(524,1705)(1499,1930)(749,2680)(524,1705)
\path(2324,3505)(1499,3130)(2549,2455)(2324,3505)
\path(2549,2455)(1499,1930)(2474,1105)(1349,805)
\path(2549,2455)(2474,1105)(4124,1105)
	(2549,2455)(3674,2830)(4124,1105)
\path(2324,3505)(2924,3355)(3674,2830)
\path(2924,3355)(2549,2455)
\end{picture}
}
\caption{\label{ideal} A piece of ideal circle pattern, with its incidence
  graph.}
\end{center}
\end{figure}

Given two of the circles in an ideal circle pattern, we consider the
angle between them only when the corresponding vertices of the incidence graph
share an edge. Thus the intersection angles between the circles correspond to
a function from the edges of the incidence graph to $(0,\pi)$.
The angle will always be measured in the complement of one disk
in the other, i.e. it will be equal to $\pi$ minus the
intersection angle measured in the intersection of the two disks. By
extension, the intersection angle between the boundaries of two disjoint but
tangent open disks is equal to $\pi$.

The notion of ideal circle pattern is not restricted to the sphere; one can
consider it in an Euclidean or a hyperbolic surface. Actually the most natural
setting is a surface with a $\C P^1$-structure, since the notion of circle
is then well defined, as well as the notion of angle. A surface with a
spherical, Euclidean or hyperbolic metric also has a canonical $\C
P^1$-structure coming from the Poincar\'e uniformization theorem. 

\pg{Ideal polyhedra.}

An {\it ideal hyperbolic polyhedron} is a convex polyhedron in $H^3$, of
finite volume, with all its vertices on the sphere at infinity. Another
possible definition is as the convex hull of a finite set of points (not all
contained in a plane) on the sphere at infinity in $H^3$. 

The possible combinatorics and dihedral angles of
ideal polyhedra were described by Andreev \cite{andreev-ideal} and Rivin
\cite{rivin-annals}. Thurston \cite{thurston-notes} realized that there was a
deep connection between ideal polyhedra and circle patterns: given an ideal
polyhedron, the boundaries of the hyperbolic planes containing its faces is an
ideal circle pattern; its incidence graph is combinatorially dual
to the 
combinatorics of the ideal polyhedron, and the exterior dihedral angles
between the faces are equal to the intersection angles between the circles. 
Using this, Thurston gave a simple
proof of the Koebe circle packing theorem based on Andreev's theorem. Thus
Rivin's theorem on ideal polyhedra can be used to describe the possible
combinatorics and dihedral angles of ideal circle patterns, in a way which is
similar to Theorem \ref{tm:ball} below. Many interesting extensions and
important references can be found in
\cite{bobenko-springborn,rivin-combi,leibon1}, in 
particular regarding the extension to higher genus surfaces.

\pg{Hyperideal circle patterns.} 

We now introduce another notion of circle patterns. Just as the ideal circle
patterns are related to ideal hyperbolic polyhedra, those other circle
patterns are related to strictly hyperideal polyhedra, so we call them {\it
  strictly hyperideal circle patterns}. 

\begin{df} \label{df:hyper}
Let $C_1, \cdots, C_N$ be a circle pattern in $S^2$, with interstices
$I_1,\cdots, I_M$. It is {\bf strictly  hyperideal} if:
\begin{itemize}
\item Each interstice has non-empty interior. 
\item For each $j\in \{ 1, \cdots, M\}$, there is an oriented circle $C'_j$,
  containing $I_j$, which is orthogonal to all the circles $C_i$ adjacent to
  $I_j$.
\item For all $i\in \{1,\cdots, N\}$ and all $j\in \{ 1,\cdots, M\}$, if $C_i$
  is not adjacent to $I_j$, then either the interior of 
  $C_i$ is disjoint from the interior of
  $C'_j$, or $C_i$ intersects $C'_j$ and their intersection angle is strictly
  larger than $\pi/2$. 
\item If $D$ is an open disk in $S^2$ such that:
  \begin{enumerate}
  \item For each $j\in \{ 1,\cdots, M\}$, either $D$ is disjoint from the
  interior of $C'_j$, or $\dr D$ has an intersection angle at least $\pi/2$ with
  $C'_j$. 
  \item $\dr D$ is orthogonal to at least $3$ of the $C'_j$.
  \end{enumerate}
then $\dr D$ is one of the $C_i$.
\end{itemize}
\end{df}

As for ideal circle patterns, this notion is not restricted to the sphere, but
can be considered for any surface with a $\C P^1$-structure.

Taking the limit of such strictly 
hyperideal circle patterns as the radii of the $C'_j$
goes to zero yields ideal circle patterns. On the other hand, taking the
limit as the intersection angles between the $C_i$ goes to $\pi$ yields a
circle packing as in the Koebe theorem quoted above.
There is a natural notion which
generalizes both the ideal and the strictly hyperideal circle patterns, and
the name ``hyperideal circle pattern'' should be kept for this more general
notion.  

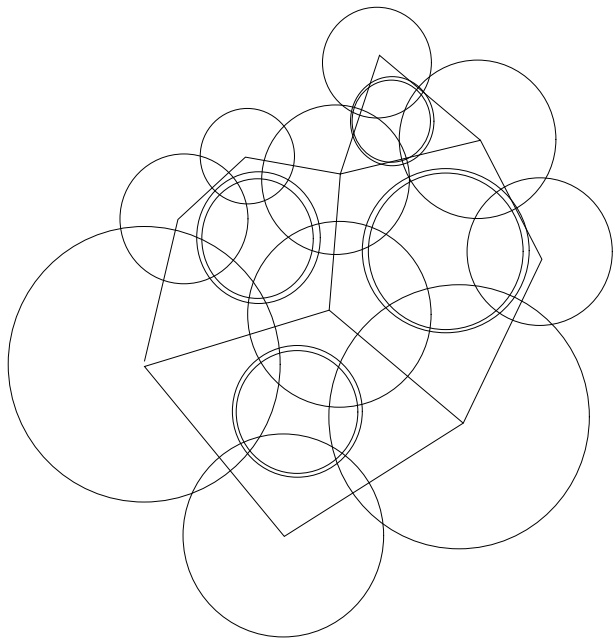
\begin{figure}[htbp]
\begin{center}
\setlength{\unitlength}{0.00083333in}
\begingroup\makeatletter\ifx\SetFigFont\undefined%
\gdef\SetFigFont#1#2#3#4#5{%
  \reset@font\fontsize{#1}{#2pt}%
  \fontfamily{#3}\fontseries{#4}\fontshape{#5}%
  \selectfont}%
\fi\endgroup%
{\renewcommand{\dashlinestretch}{30}
\begin{picture}(3786,3959)(0,-10)
\put(2076,2020){\ellipse{1146}{1162}}
\put(2051,2861){\ellipse{926}{938}}
\put(2822,1381){\ellipse{1628}{1650}}
\put(858,1709){\ellipse{1700}{1722}}
\put(1502,3008){\ellipse{590}{598}}
\put(1106,2617){\ellipse{798}{810}}
\put(2938,3113){\ellipse{978}{992}}
\put(3324,2412){\ellipse{908}{922}}
\put(1726,641){\ellipse{1252}{1268}}
\put(2309,3592){\ellipse{682}{692}}
\put(2739,2415){\ellipse{968}{980}}
\put(1811,1411){\ellipse{760}{770}}
\put(1813,1414){\ellipse{812}{824}}
\put(2739,2418){\ellipse{1042}{1024}}
\put(2405,3226){\ellipse{522}{558}}
\put(1571,2499){\ellipse{772}{824}}
\put(1565,2493){\ellipse{698}{748}}
\put(2399,3227){\ellipse{488}{514}}
\path(861,1729)(1070,2614)(1490,3002)
	(2082,2896)(2326,3638)(2955,3108)
	(3339,2365)(2850,1341)(1733,634)
	(861,1695)(2013,2048)(2082,2896)(2955,3108)
\path(2013,2048)(2850,1341)
\end{picture}
}
\caption{\label{hyper} A piece of hyperideal circle pattern, with its incidence
  graph. The circles corresponding to the interstices are  doubled. The angles
  between the circles of the two families should be $\pi/2$.} 
\end{center}
\end{figure}

Given a strictly hyperideal circle pattern, its {\it incidence graph} is the
1-skeleton of the cellular decomposition of $S^2$
defined similarly as for ideal circle patterns, which has:
\begin{itemize}
\item One face for each of the $I_j, 1\leq j\leq M$.
\item One vertex for each of the $C_i, 1\leq i\leq N$.
\item One edge between two faces, corresponding to $I_k$ and $I_l$, 
whenever there are two circles $C_i$ and $C_j$ which are both orthogonal to
$C'_k$ and $C'_l$.
\end{itemize}
We will consider the intersection angles between the circles $C_i$, defined as
explained above for ideal circle patterns. 

\pg{Some examples.}

Before stating the main result, we give three simpler examples, two of them
corresponding to already well understood cases. The first is a direct
consequence of a recent result of Bao and Bonahon \cite{bao-bonahon}.

\begin{thm} \label{tm:ball}
Let $\Gamma$ be the 1-skeleton of a polytopal cellular decomposition of
$S^2$, and let $w:\Gamma_1\rightarrow (0,\pi)$ be a map on the set of edges of
$\Gamma$. There exists a strictly hyperideal circle pattern on $S^2$ with
incidence graph $\Gamma$ and intersection angles given by $w$ if and only if:
\begin{enumerate}
\item For each simple closed curve $\gamma$ in $\Gamma$,
  the sum of the values of 
  $w$ on the edges of $\gamma$ is strictly larger than $2\pi$.
\item For each open path $\gamma$ in $\Gamma$, which begins and
  ends on the boundary of a face $f$ but is not contained in $f$, 
  the sum of the values of 
  $w$ on the edges of $\gamma$ is strictly larger than $\pi$.
\end{enumerate}
This strictly hyperideal circle pattern is then unique, up to the M\"obius
transformations of $S^2$ .
\end{thm}

There is a similar result in the case of the torus. 

\begin{thm} \label{tm:torus}
Let $\Gamma$ be the 1-skeleton of a cellular decomposition of
$T^2$, and let $w:\Gamma_1\rightarrow (0,\pi)$ be a map on the set of edges of
$\Gamma$. There exists a flat metric $g_0$ on $T^2$ and a strictly hyperideal
circle pattern $C$ on $(T^2, g_0)$ with
incidence graph $\Gamma$ and intersection angles given by $w$ if and only if:
\begin{enumerate}
\item For each simple, homotopically trivial 
  closed path $\gamma$ in $\Gamma$, the sum of the values of 
  $w$ on the edges of $\gamma$ is strictly larger than $2\pi$.
\item For each open path $\gamma$ in $\Gamma$, which begins and
  ends on the boundary of a face $f$, is homotopic to a segment in that face,
  but is not contained in $f$, 
  the sum of the values of $w$ on the edges of $\gamma$ is strictly larger than
  $\pi$. 
\end{enumerate}
Then $(g_0, C)$ is unique, up to the homotheties of $g_0$.
\end{thm}

The same statement holds on surfaces of higher genus, it is a direct
consequence of (a special case of) a recent result of Rousset \cite{rousset1}.

\begin{thm} \label{tm:higher}
Let $S_g$ be a closed surface of genus $g\geq 2$, let 
$\Gamma$ be the 1-skeleton of a cellular decomposition of
$S_g$, and let $w:\Gamma_1\rightarrow (0,\pi)$ be a map on the set of edges of
$\Gamma$. There exists a hyperbolic metric $g_0$ on $S_g$ and a strictly
hyperideal circle pattern $C$ on $(S_g, g_0)$ with 
incidence graph $\Gamma$ and intersection angles given by $w$ if and only if:
\begin{enumerate}
\item For each simple, homotopically trivial 
closed path $\gamma$ in $\Gamma$,
the sum of the values of 
  $w$ on the edges of $\gamma$ is strictly larger than $2\pi$.
\item For each open path $\gamma$ in $\Gamma$, which begins and
  ends on the boundary of a face $f$, is homotopic to a segment $f$, but is
  not contained in $f$, 
  the sum of the values of $w$ on the edges of $\gamma$ is strictly larger than
  $\pi$. 
\end{enumerate}
Then $(g_0, C)$ is unique.
\end{thm}

\pg{The main result.} 

The 3 theorems quoted above are basically simple consequences of the following
more complex but much more general statement. 
We now consider a compact 3-manifold with non-empty 
boundary $M$, whose interior admits a
complete hyperbolic metric; $M$, like all the
manifolds that we consider in this paper, will be oriented. According to a
result of Thurston \cite{thurston-bulletin}, the existence of a complete
hyperbolic metric on the interior
of $M$ is equivalent to a simple topological condition: 
that $M$ is irreducible and homotopically atoroidal, and
is not the interval bundle over the Klein bottle.

We call $\dr'M$ the union of the connected components
of $\dr M$ which are not tori, except if $M$ is a solid torus --- then
$\dr'M=\dr M$ --- or if $M$ is the product of a torus by an interval --- then
$\dr 'M$ is one connected component of $\dr M$.

In this setting, the complete hyperbolic metrics on $M$ are 
{\it geometrically finite}; they contain a finite volume subset $C$ which 
is {\it convex} in the (strong) sense that, given {\it any} 
geodesic segment $\gamma$ with
endpoints in $C$, $\gamma$ is contained in $C$ (see
\cite{thurston-notes}). For each such metric,
$M$  is isometric to the quotient of $H^3$ by a discrete group $G$ acting by
isometries. Then $\dr'M$ is the quotient of the discontinuity domain of $G$ by
$G$. Since $G$ acts on $S^2$ by M\"obius transformations, $\dr'M$ has a
natural $\C P^1$-structures. Note that, in general, 
among the $\C
P^1$-structures on $\dr 'M$, only a small subset are obtained from a
geometrically finite hyperbolic metric on $M$.

Finally, before stating the result, we need a small restriction on the kind of
cellular decomposition which can be realized.

\begin{df}
Let $M$ be a compact 3-manifold with boundary, whose interior admits a
complete hyperbolic metric. Let $\Sigma$ be a cellular
decomposition of $\dr'M$. $\Sigma$ is {\it
  proper} if there is no essential disk $D$ in $M$ such that $\dr D\subset
\dr'M$ and that $\dr D$ intersects the closure of at most two cells.  
\end{df}

\begin{thm} \label{tm:circles}
Let $M$ be a compact (orientable) 
3-manifold with boundary whose interior admits a
complete hyperbolic metric. 
Let $\Gamma$ be the 1-skeleton of a proper cellular decomposition $\Sigma$ of
$\dr'M$, and let 
$w:\Gamma_1\rightarrow (0,\pi)$ be a map on the set of edges of 
$\Gamma$. There exists a couple $(\sigma, C)$, where $\sigma$ is a $\C
P^1$-structure on $\dr'M$ induced by a geometrically finite hyperbolic metric
on $M$ and $C$ is a strictly 
hyperideal circle pattern on $(\dr'M, \sigma)$ with 
incidence graph $\Gamma$ and intersection angles given by $w$, if and only if:
\begin{enumerate}
\item The cellular decomposition $\Sigma$ is proper.
\item For each simple closed path $\gamma$ in $\Gamma$,
  homotopically trivial in $M$, the sum of the values of 
  $w$ on the edges of $\gamma$ is strictly larger than $2\pi$.
\item For each open path $\gamma$ in $\Gamma$, which begins and
  ends on the boundary of a face $f$, is homotopic in $M$ 
  to a segment in $f$, but is not contained in $f$,
  the sum of the values of $w$ on the edges of
  $\gamma$ is strictly larger than $\pi$. 
\end{enumerate}
Then $(\sigma, C)$ is unique.
\end{thm}

Theorems \ref{tm:ball}, \ref{tm:torus} and \ref{tm:higher}
are direct consequences of Theorem \ref{tm:circles}. Theorem \ref{tm:ball} is
obtained when $M$ is a ball, and Theorem \ref{tm:torus} when $M$ is the
product of a torus by $\R$. For Theorem
\ref{tm:higher}, one should take as $M$ the product of a surface of genus at
least $2$ by an interval, with a graph $\Gamma$ and a function $w$ which are
the same on both connected components of the boundary; the uniqueness in
Theorem \ref{tm:circles} then shows that the $\C P^1$-structure obtained on
each connected component of the boundary is the same, so that it is the $\C
P^1$-structure induced by a hyperbolic metric on the surface considered.

There are other simple consequences of Theorem \ref{tm:circles} beyond the
three statements given above. For instance, one can consider two cellular
decompositions of a surface $\Sigma$ 
of genus at least $2$, along with some angles on
the edges, and obtain from Theorem \ref{tm:circles} a unique quasi-fuchsian
metric on $\Sigma\times \R$,
along with one strictly hyperideal circle pattern on each boundary
component. 

\pg{Relation to other results.}

Theorem \ref{tm:circles} is similar to the main result of \cite{hphm},
although \cite{hphm} is about hyperideal polyhedra rather than hyperideal
circle patterns. Theorem is \ref{tm:circles}
less general than the main result of \cite{hphm} insofar as 
the circle patterns considered are strictly hyperideal
rather than hyperideal (a more general notion), but also more general since
$M$ can have compressible boundary and $\dr M$ is allowed to have toric
components. The proof of Theorem \ref{tm:circles}, however, is completely
different from the one given in \cite{hphm}; here we show that Theorem
\ref{tm:circles} is a direct consequence of a result of Otal
\cite{otal,bonahon-otal}, 
which itself uses crucially a result of Hodgson and
Kerckhoff \cite{HK}. In \cite{hphm}, the proof is direct and based, among
other things, on some properties of the volume of hyperideal polyhedra. The
method of \cite{hphm} yields additional informations, in particular concerning
the ``induced metrics'' on the boundaries of the manifolds with hyperideal
boundary that appear here in section 2.

\section{Circle patterns and hyperideal polyhedra}

In this section we define ``manifolds with strictly hyperideal boundary''
and state a theorem, concerning them, which is basically equivalent to Theorem
\ref{tm:circles}. This theorem is proved in the next section, using
\cite{otal,bonahon-otal}. 

\pg{Hyperideal polyhedra.}

Recall that the Klein model of $H^3$ is a map $\phi:H^3\rightarrow B^3$, where
$B^3$ is the open unit ball in $\R^3$, which is projective, i.e. it sends the
hyperbolic geodesics to the segments. A {\it hyperideal polyhedron} is a
hyperbolic polyhedron which is the inverse image, in $H^3$, of a polyhedron
$P\subset \R^3$ with all its vertices outside $B^3$ but all its edges
intersecting $B^3$. It is {\it strictly hyperideal} if $P$ has no vertex
on the boundary of $B^3$.

Given a point $v\in \R^3\setminus \overline{B^3}$, its {\it polar dual} is the
plane, noted $v^*$, which contains 
the set of points $x\in S^2$ such that the line going
through $x$ and $v$ is tangent to $S^2$. (It can also be defined using the
polarity with respect to a bilinear form of signature $(3,1)$ on $\R^4$, which
explains the terminology.) An important property is that, for any $x\in
v^*\cap B^3$, the intersection with $B^3$ of the line going through $x$ and
$v$, considered as a hyperbolic geodesic, is orthogonal to $v^*\cap B^3$,
considered as a hyperbolic plane. Conversely, for any plane $P\subset \R^3$
intersecting $B^3$, there is
a unique point $P^*\in \R^3\setminus B^3$ such that $P$ is the plane dual to
$v$. This notion of duality has many interesting
application, it is related to the hyperbolic-de Sitter duality used in
particular by Rivin \cite{Ri} and Rivin and Hodgson \cite{RH} to obtain
beautiful results on compact hyperbolic polyhedra. 
To simplify statements a little, we will sometimes talk about the point which
is dual to an oriented circle in $S^2$; it is the point dual to the plane
which contains the circle. 

Using this notion, we can reformulate the definition of a hyperideal polyhedron
without reference to the projective model of $H^3$; a hyperideal polyhedron is
a convex hyperbolic polyhedron (i.e. the intersection of a finite number of
half-planes) with no vertex such that, for each end,
either the end has finite volume, or there exists a hyperbolic plane which is
orthogonal to all the edges going to infinity in it. It is strictly hyperideal
if all ends have infinite volume. 

Moreover, given a hyperideal polyhedron $P$, 
one can {\it truncate} it (following \cite{bao-bonahon}) 
by cutting off each end of infinite volume
by the plane which is dual to the corresponding vertex; one obtains
in this way a finite volume polyhedron $P_t$ 
(which is compact if $P$ is strictly
hyperideal). Its faces are either faces of $P$, or faces which are orthogonal
to the adjacent faces.

\pg{Ideal polyhedra and ideal circle patterns.}

We first recall the correspondence between ideal hyperbolic polyhedra and
ideal circle patterns. It is a direct consequence of the following
proposition, which is classical, see e.g. \cite{handbook}.

\begin{prop} \label{pr:ideal}
Let $x_1, \cdots, x_M\in S^2$ be distinct points, $M\geq 3$. 
There exists a unique family of oriented
circles $C_1, \cdots, C_N$ such that:
\begin{itemize}
\item The $C_i$ bound closed disks which cover $S^2$.
\item None of the $C_i$ bounds an open disk containing any of the $x_j$.
\item Each of the $C_i$ contains at least 3 of the $x_j$.
\item If $D$ is an open disk which contains none of the $x_j$ but such that
  $\dr D$ contains at least $3$ of them, then $D$ is the interior of one of
  the $C_i$.
\end{itemize}
The $C_j$ are obtained as the intersections with $S^2$ of the planes containing
the faces of the convex hull, in $\R^3$, of the points $x_1, \cdots, x_M$.
\end{prop}

This proposition explains the relationship between ideal polyhedra and ideal
circle patterns. Given a (convex) ideal polyhedron, it is clear that the set
of oriented circles associated to its faces is an ideal circle
pattern. Conversely, given an ideal circle pattern $C_1, \cdots, C_N$, let
$x_1, \cdots, x_M$ be its interstices. 
Then Definition \ref{df:ideal} implies that $C_1, \cdots, C_N$ is
the set of oriented circles associated by Proposition
\ref{pr:ideal} to $x_1, \cdots, x_M$. Therefore $C_1, \cdots, C_N$ are the
boundaries of the planes containing the faces of the ideal polyhedron which is
the convex hull of the $x_j$.

\pg{Hyperideal polyhedra and hyperideal circle patterns.}

The same procedure can be applied for strictly hyperideal circle patterns,
based on the following analog of Proposition \ref{pr:ideal} (see
\cite{handbook}). 

\begin{prop} \label{pr:hyper}
Let $C'_1, \cdots, C'_M$ be a family of oriented circles in $S^2$, 
bounding disjoint
closed disks ($M\geq 3$), each disk being strictly smaller than a hemisphere. 
There exists a unique family of oriented circles $C_1, \cdots, C_N$ such that:
\begin{itemize}
\item The closed disks bounded by the $C_i$ cover the complement of the closed
  disks bounded by the $C'_j$.
\item None of the $C_i$ contains any of the $C'_j$ in its interior, or makes
  an angle strictly less than $\pi/2$ with any of the $C'_j$.
\item Each of the $C_i$ is orthogonal to at least $3$ of the $C'_j$.
\item If $D$ is an open disk such that:
  \begin{enumerate}
  \item For each $j\in \{ 1, \cdots, M\}$, either $D$ is disjoint from the
  interior of $C'_j$, or $\dr D$ makes an angle at least $\pi/2$ with $C'_j$.
  \item $\dr D$ is orthogonal to at least $3$ of the $C'_j$.
  \end{enumerate}
then $D$ is the interior of one of the $C_i$.
\end{itemize}
The $C_i$ are obtained as the intersections with $S^2$ of the planes containing
the faces of the convex hull, in $\R^3$, of the points $x_1, \cdots, x_M$
which are dual to the circles $C'_1, \cdots, C'_M$.
\end{prop}

Clearly the hypothesis that the circles $C'_j$ bound disks smaller than
hemispheres is not crucial, it appears here in the convex hull construction 
because we do things in the
Euclidean space rather than in the sphere, where statements would be a little
more complicated. Note also that there is a common generalization of
Proposition \ref{pr:ideal} and Proposition \ref{pr:hyper}, see
\cite{handbook}. 

Given a  strictly hyperideal polyhedron $P$, 
one can consider the oriented planes
containing its faces, and then their boundaries, which are oriented circles in
$S^2$. Those circles clearly make up a strictly hyperideal circle pattern,
with the connected components of the 
complement of the disks corresponding to the vertices of
$P$; the circles $C'_j$ which appear in the definition of a strictly 
hyperideal circle
pattern are the boundaries of the hyperbolic planes dual to the strictly
hyperideal vertices. Moreover, the incidence graph of this
circle pattern is dual to the combinatorics as the polyhedron $P$, because two
faces of the incidence graph --- corresponding to two connected components
of the complement of the disks --- are adjacent if and only if the
corresponding circles are both orthogonal to two circles of the pattern,
i.e. if and only if the corresponding hyperideal vertices of $P$ are adjacent.

Conversely, let $C=(C_1, \cdots, C_N)$ be a strictly hyperideal circle
pattern. Let $C'_1, \cdots, C'_M$ be the oriented circles corresponding to the
interstices. 
Since those circles bound disjoint disks, we can apply a M\"obius
transformation to ensure that all those disks are strictly smaller than
hemispheres. 

Definition \ref{df:hyper} shows that the circles $C_1, \cdots, C_N$ are
precisely the circles associated to $C'_1, \cdots, C'_M$ by Proposition
\ref{pr:hyper}. Therefore, the $C_i$ are the boundaries of the planes
containing the faces of a strictly 
hyperideal polyhedron, which is the convex hull of
the points dual to the $C'_j$.
 
\pg{Manifolds with hyperideal boundary.}

It is necessary below to consider objects, more general than
hyperideal polyhedra, which are basically submanifolds $M$ of a geometrically
finite hyperbolic 3-manifold $N$ such that the boundary of $M$ in $N$ is
locally like a strictly hyperideal polyhedron in $H^3$.

A {\it geometrically finite} hyperbolic 3-manifold is the interior of 
a compact 3-manifold
with boundary $N$ with a complete hyperbolic metric, which contains a non-empty
submanifold $M\subset N$, of finite volume, which is {\it convex} in the
(strong) sense that any
geodesic segment in $N$ with endpoints in $M$ actually remains in $M$. 
Moreover, if $N$ is not a ball or the product of a torus by an interval, 
there exists a smallest such subset $M_0$ (smallest
with respect to the inclusion) and $M_0$ is called the {\it convex core} of
$N$. If $N$ is not a ball, a solid torus or the product of the torus by an
interval,
then the boundary of $M_0$ in $N$ is homeomorphic to the union of the connected
components of $\dr N$ which are not tori, while the tori in the boundary of
$N$ correspond to cusps in $M_0$ (see e.g. 
\cite{thurston-bulletin,thurston-notes} for much more on this).

\begin{df}
A {\it hyperbolic manifold with strictly
hyperideal boundary} is the interior of a compact 3-manifold
with boundary $M$ with a hyperbolic metric, such that there exists an
isometric embedding $\phi:M\rightarrow N$, where $N$ is a geometrically finite
hyperbolic 3-manifold, with the property that:
\begin{itemize}
\item $\phi(M)$ is convex in $N$ in the strong sense defined above: any
  geodesic segment of $N$ with endpoints in $\phi(M)$ actually remains in
  $\phi(M)$.  
\item For any embedded hyperbolic ball $B\subset N$, there exists a strictly
  hyperideal polyhedron $P\subset H^3$ and a ball $B'\subset H^3$ such that
  the interior of $B'\cap P$ is isometric to $B\cap \phi(M)$.
\item The boundary of $\phi(M)$ in $N$ contains no closed geodesic of $N$.
\end{itemize}
Given $M$, $N$ is uniquely defined and will be called the {\it extension} of
$M$. 
\end{df}

In many cases it is more convenient to consider the closure of $M$ in $N$,
rather than $M$ as defined above.
The condition on the closed geodesics is equivalent to the fact
that the boundary of $M$ in $N$
does not intersect the convex core of $N$. 
There is another possible definition, using the notion of hyperideal point in
a geometrically finite manifold. A simple example of a strictly hyperideal
manifold is that, when $N$ is the
ball, $M$ is simply a strictly hyperideal polyhedron.

\pg{Dihedral angles of manifolds with hyperideal boundary.}

Consider a hyperbolic manifold $M$ with strictly hyperideal boundary. 
By definition, its boundary $\dr M$ is a ``polyhedral'' surface in the
extension $N$ of $M$. It has a finite number of faces and a finite number of
edges. Each face is isometric to the interior of a strictly hyperideal polygon
in $H^2$, while the edges are complete geodesics. The boundary combinatorics
of $M$ determines (combinatorially) a cellular decomposition of $\dr'N$ ---
with one face for each face of $\dr M$, and one vertex for each strictly
hyperideal vertex of $M$ (or, equivalently, for each end of $\dr M$). 

For each edge $e$ of $\dr M$, we define the {\it exterior dihedral angle} of
$\dr M$ at $e$ to be $\pi$ minus the angle, measured in $M$, between the two
faces containing $e$. Another possible definition is as the angle, measured on
$e$, between the exterior normal vectors to the two faces of $\dr M$
containing $e$. The possible combinatorics and boundary angles of the strictly
hyperideal manifolds which have as their extension a given geometrically
finite hyperbolic 3-manifold $N$ are described in the following statement. It
has non-empty intersection with the results of
\cite{bao-bonahon,rousset1,hphm}.

\begin{thm} \label{tm:poly}
Let $N$ be compact 3-manifold with boundary, whose interior admits a complete
hyperbolic metric. Let $\Gamma$ be a
graph, embedded in $\dr'N$, which is the 1-skeleton of a cellular
decomposition $\Sigma$ of $\dr'N$, and
let $w:\Gamma_1\rightarrow (0,\pi)$ be a map from
the set of edges of $\Gamma$ to $(0,\pi)$. There exists a strictly hyperideal
hyperbolic manifold $M$, which can be isometrically embedded in $N$ with some
geometrically finite metric, with boundary
combinatorics given by $\Gamma$ and exterior dihedral angles described by $w$,
if and only if:
\begin{enumerate}
\item The cellular decomposition $\Sigma$ is proper.
\item For each simple 
  closed path $\gamma$ in the dual graph $\Gamma^*$ of $\Gamma$,
  which is homotopically trivial in $M$, the sum of the values of $w$ on the
  edges of $\gamma$ is strictly larger than $2\pi$.
\item For each open path $c$ in $\Gamma^*$ which begins and ends on the same
  face $f$ of $\Gamma^*$, is not contained in the boundary of $f$, but is
  homotopic in $M$ to a segment in $f$, 
  the sum of the values of $w$ on the edges of $c$ is strictly larger than
  $\pi$. 
\end{enumerate}
This manifold $M$ is then unique. 
\end{thm}

Note that the existence of a function $w$ satisfying the requested conditions
implies constraints on $\Gamma$; for instance, if $M$ is a  ball, $\Gamma$ has
to be the 1-skeleton of a polytopal cellular decomposition (see
\cite{bobenko-springborn,hphm} for more details on this).

\pg{Proof of Theorem \ref{tm:circles} from Theorem \ref{tm:poly}.}

Let $N$, $\Gamma$ and $w$ be given as in the statement of Theorem
\ref{tm:circles}. Apply Theorem \ref{tm:poly} to obtain a strictly hyperideal
hyperbolic metric $h$ on $N$ and a strictly hyperideal manifold $M$
isometrically embedded in $(N,h)$. 
Let $\sigma$ be the $\C P^1$-structure induced on
$\dr'N$ by $(N,h)$, and let $C$ be the strictly
hyperideal circle pattern defined on $\dr'N$ by the boundaries of the planes
containing the faces of $M$. By construction, $C$ has the combinatorics
dual to the combinatorics of $\Gamma$ and the intersection angles determined by
$w$. Moreover, the uniqueness in Theorem \ref{tm:circles} will follow from the
uniqueness in Theorem \ref{tm:poly} if we show that any strictly hyperideal
circle pattern in $\dr'N$, for a $\C P^1$-structures which is induced by a
geometrically finite metric on $N$, comes from a strictly hyperideal
hyperbolic manifold isometrically embedded in $N$ (with some complete
hyperbolic metric, not necessarily $h$).

Let $\sigma$ be a $\C P^1$-structure on $\dr'N$ induced by a geometrically
finite hyperbolic metric $g$ on $N$, and let $C$ be a strictly  hyperideal
circle pattern in $\dr'N$ for $\sigma$. Let $\gt$ be the lift of $g$ to the
universal cover $\Nt$ of $N$; then $(\Nt,\gt)$ is isometrically identified
with $H^3$, and the fundamental group $\pi_1N$ of $N$ acts on $H^3$ by
isometries. Let $\Lambda\subset S^2$ be its limit set. 

Under the projection from $\Nt$ to $N$, $\dr'N$ lifts to the complement of
$\Lambda$ in $S^2$. The strictly hyperideal circle pattern $C$ lifts to an
infinite set $\Ct$ of circles in $S^2\setminus \Lambda$, which accumulates
close to $\Lambda$ --- each neighborhood of a point of $\Lambda$ contains an
infinite number of circles. However $\Ct$ is still a strictly hyperideal
circle pattern in the sense of Definition \ref{df:hyper} (except of course
for the finiteness condition). 

It is well known that Proposition \ref{pr:hyper} extends to the
context where the family of disjoint circles is infinite, see
\cite{handbook} (it can also be proved by a simple limiting argument based on
Proposition \ref{pr:hyper}). It follows, as
in the case of circle patterns on the sphere, that $\Ct$ is the set of
boundaries of the planes containing the faces of the convex hull of an infinite
set of points $S\subset \R^3$, which accumulates near $\Lambda$.

By the uniqueness in the construction, $S$ is invariant under the action of
$\pi_1N$ on $\R^3$ (which is by projective transformations preserving
$S^2$). Taking the quotient of the convex hull of $S$ by the action of 
$\pi_1N$ defines, in $N$, a manifold with strictly
hyperideal boundary, whose combinatorics and boundary angles are as needed. 
Thus Theorem \ref{tm:circles} follows from Theorem \ref{tm:poly}.

\section{Hyperideal polyhedra and hyperbolic convex cores}

It remains, in this section, to prove Theorem \ref{tm:poly}. We will show here
that it is a simple consequence of a result of Otal
\cite{otal,bonahon-otal}.

\pg{Convex cores of hyperbolic manifolds.}

Let $N$ be a geometrically finite hyperbolic 3-manifold (not a ball, a solid
torus or the
product of a torus by an interval). We have defined above
the convex core of $N$ as the smallest non-empty 
subset of $N$ which is {\it convex} in the sense that, given a geodesic
segment $\gamma$ in $N$ with endpoints in $C$, $\gamma\subset C$. The convex
core of $N$ has finite volume (by definition of a geometrically finite
hyperbolic manifold). 

The boundary of the convex core of $N$ is a non-smooth
surface in $N$ which is locally
convex; it has a hyperbolic induced metric, but is ``bent'' along a
lamination (see \cite{thurston-notes}), called the ``pleating
lamination''. The precise way the bending occurs is described by a transverse
measure defined on the pleating lamination. We
are interested here only in a very special case, where the support of the
pleating lamination is a disjoint union of simple closed curves; in this case
we say that the pleating lamination is {\it rational}.

In this case, the boundary of the convex core is a polyhedral-like surface. It
has totally geodesic faces, edges which are closed geodesics (corresponding to
the curves in the support of the pleating lamination) and no vertex. To each
edge is associated a number --- describing the weight of the corresponding
leaf for the transverse measure of the pleating lamination --- which is simply
the angle between the oriented normals of the ``faces'' on the two sides; see
e.g. \cite{epstein-marden,bonahon-geodesic} for more on the local geometry of
the boundary of the convex core. 

It follows from its definition that the convex core of $N$ has the same
topology as $N$.
Since the pleating lamination and its bending measure are topological in
nature, we can consider them on the boundary of $N$ rather than the boundary
of its convex core. 

\pg{A theorem of Otal.}

The possible transverse measures of the rational pleating laminations of
hyperbolic convex cores are precisely described by the following result. 

\begin{thm}[Otal \cite{otal,bonahon-otal}] \label{tm:otal}
Let $N$ be a compact 3-manifold with boundary, whose interior admits a
complete hyperbolic metric. Let $\alpha$ be a measured lamination on $\dr N$
whose leaves are closed curves. There exists a non-fuchsian geometrically
finite metric $g$ such that $\alpha$ is the measured pleating lamination of
the boundary of the convex core if and only if the following conditions are
satisfied:
\begin{enumerate}
\item Each leaf of $\alpha$ has weight at most $\pi$.
\item For all essential annulus or M\"obius strip $A$ in $N$, the total weight
  of $\dr A$ is strictly positive.
\item For all essential disk $D$ in $N$, the total weight of $\dr D$ is
  strictly larger than $2\pi$.
\end{enumerate}
This metric $g$ is then unique (up to isotopy in $N$).
\end{thm}
 
The proof of this difficult result uses among other things a result of Hodgson
and Kerckhoff \cite{HK} on the rigidity of hyperbolic cone-manifolds. 

\pg{Proof of Theorem \ref{tm:poly} from Theorem \ref{tm:otal}.}

First note that condition (1) of Theorem \ref{tm:poly}, that the cellular
decomposition of $\dr'N$ is proper, is necessary. If it were not satisfied, and
if $M$ were a manifold with hyperideal boundary as described in the statement
of Theorem \ref{tm:poly}, the boundary of the universal cover $\Mt$ of
$M$, considered as a complete convex polyhedral surface in $H^3$, 
would contain a
non-trivial closed path contained in the union of the closures of two faces. 
This is clearly impossible by convexity.  

The rest of the proof proceeds in several steps. 

\medskip
\noindent {\bf 1. The uniqueness in Theorem \ref{tm:poly}.}
Let $N, \Gamma$ and $w$ be as in the statement of Theorem \ref{tm:poly}.
Suppose first that, as stated in the theorem, there exists a strictly 
hyperideal hyperbolic metric $g$ on $N$ such that the boundary combinatorics
of $N$ is given by $\Gamma$ and its dihedral angles by $w$. 
Let $N_t$ be the truncated version of $N$. So $N_t$ is a finite volume
hyperbolic manifold
with convex, polyhedral boundary; its boundary has one face for each face of
$\Gamma$, and one --- induced by the truncation --- for each vertex of
$\Gamma$ (i.e. for each hyperideal vertex of $N$). 
Moreover, each face corresponding to a vertex of $\Gamma$ is orthogonal to
each face corresponding to an adjacent face of $\Gamma$.

Consider another copy of $N_t$, which we call $N'_t$, and glue $N_t$ to $N'_t$
along all the couples of faces --- one on $N_t$ and the other on $N'_t$ ---
corresponding to the vertices of $\Gamma$; let $C$ be the resulting
hyperbolic manifold. By construction, $C$ has finite volume --- since $N_t$
and $N'_t$ have finite volume --- and it has a ``polyhedral'' boundary, with
no vertex. 

Each face $f$ of $\Gamma$ has two corresponding faces, one, say $F$, on $\dr
N_t$, and the other, say $F'$, on $\dr N'_t$. But, for each vertex $v$ of
$\Gamma$ adjacent to $f$, the faces on $\dr N_t$ and $\dr N'_t$
corresponding to $v$ are orthogonal to $F$ and $F'$, respectively. If follows
that $F$ and $F'$ are two parts of the same face of $C$. Each edge of $F$
with, on the other side, a face of $\dr N_t$ corresponding to a vertex of
$\Gamma$ ``disappears'' in $C$, while each edge of $F$ with on the other side
a face of $\dr N_t$ contained in a face of $\dr N$ becomes, in $C$, half of an
edge of $\dr C$ which is a closed curve. 

The following description of $\dr C$ follows from those considerations. It
has:
\begin{itemize}
\item One face $F_c$ for each face $f$ of $\Gamma$, and $F_c$ is topologically
  a sphere with $v$ disks removed, where $v$ is the number of vertices of $f$.
\item One edge for each edge of $\Gamma$, and each of those edges is a closed
  curve. The (exterior) dihedral angle at each edge is equal to the value of
  $w$ on the corresponding edge of $\Gamma$.
\end{itemize}
It cleary follows from this description that $C$ is the convex core of a
geometrically finite hyperbolic 3-manifolds, which has a rational pleating
lamination. The uniqueness in Theorem \ref{tm:otal} shows that $C$ is uniquely
determined by $\Gamma$ and by $w$, and the uniqueness in Theorem \ref{tm:poly}
follows. 

\medskip
\noindent {\bf 2. Conditions (2) and (3) of Theorem \ref{tm:poly} are
  necessary.} 
For condition (2), consider a closed curve $\gamma$ as in the statement of
condition (2), so that $\gamma$ is a sequence of faces of $\dr M$, with two
consecutive faces sharing an edge. Associate to $\gamma$ a curve $\gamma'$ in
$\dr M$, which intersects the same faces in the same order. 
Following the construction made in step (1), $\gamma'$ can be
considered as a closed curve in the boundary of $C$. Since $\gamma$ is simple,
this curve is homotopically non-trivial in $\dr C$, while it is homotopically
trivial in $C$. But $\dr C$ lifts to a
complete, convex, polyhedral surface in $H^3$. It is known that, for any
non-trivial closed curve on such a surface, the sum of the exterior 
dihedral angles of the edges crossed by the closed curve is strictly larger
than $2\pi$, see \cite{RH,CD,these,shu}. Condition (2) follows.

For condition (3), consider an open curve $c$ as in the statement of condition
(3). In the manifold $C$ constructed in step (1) above, consider two copies of
$c$, one on $\dr N_t$, and the other on $\dr N'_t$; they can be glued at their
endpoints to obtain a closed, non-trivial curve in $C$, with boundary in $\dr
C$. The same argument as for condition (2) thus shows that
the sum of the exterior dihedral angles of the edges of $\dr C$
crossed by this closed path has to be strictly larger than $2\pi$, and the
statement of condition (3) follows.

\medskip
\noindent {\bf 3. Construction of a manifold with a rational 
measured lamination.}
To prove the existence of $M$, we consider a 3-manifold $\cC$ along with a
rational measured lamination $\alpha$ 
on $\dr \cC$ corresponding to the manifold $C$
obtained above. So $\cC$ is obtained by gluing two copies of $N$, which we
call $N$ and $N'$,
along a set of disks, one for each vertex of $\Gamma$. Since the interior of
$N$ has a complete hyperbolic metric, $N$ is irreducible and homotopically
atoroidal; since $\cC$ is obtained by gluing two copies of $N$ along some disks
in the boundaries, it follows that $\cC$ is also irreducible and homotopically
atoroidal (and $\cC$ is not the interval bundle over the Klein bottle), so
that $\cC$ admits a complete hyperbolic metric.

The lamination $\alpha$
has one closed curve for each edge of $\Gamma$, so that is has a finite number
of closed leaves. The transverse measure we consider is the one given, for
each leaf of $\alpha$, by the value of $w$ on the corresponding edge of
$\Gamma$. Clearly hypothesis (1) of Theorem \ref{tm:otal} is satisfied.

\medskip
\noindent {\bf 4. Proof that hypothesis (3) of Theorem \ref{tm:otal} is
  satisfied.}
Let $D$ be an essential disk in $\cC$. We can suppose (applying a small
deformation if necessary) that its boundary $\dr D\subset \dr \cC$ is
transverse to the edges of $\dr \cC$. In addition, it is possible to make two
kinds of deformations of $D$ to obtain a simpler picture.
\begin{itemize}
\item Suppose that some connected component $D'$ of $D\cap N$ is a (closed)
  disk in the interior of $D$. Then $\dr D'$ is contained 
  in one of the disks along which $N$ is glued to $N'$. Since $N$ is
  irreducible, $D'$ can be deformed to the disk in
  $N\cap N'$ bounded by $\dr D'$. It is then possible to ``push'' a
  neighborhood of $D'$ in $D$ into $N'$. The same argument can be applied to a
  connected component $D'$ of $D\cap N'$ under the same hypothesis.
\item Suppose that some connected component $D'$ of $D\cap N$ is such that
  $\dr D'=c\cup c'$, where $c\subset \dr D\cap \dr N$ and $c'\subset N\cap N'$
  are connected curves. Suppose moreover that $c$ remains in the union of the
  faces of $\dr N$ which are adjacent to the disk of $N\cap N'$ containing
  $c'$. Then the same argument shows that a neighborhood of $D'$ in $D$ can be
  ``pushed'' into $N'$. Again the same argument works if $D'\subset N'$. 
\end{itemize}
Finally, perturbing $D$ a little if necessary, we suppose that $D$ is
transverse to the surface $N\cap N'$ between $N$ and $N'$. 

We associate to $\dr D$ a closed path $\gamma$ in the dual graph $\Gamma^*$ of
$\Gamma$, which follows the edges dual to the edges crossed by $\dr D$. 
The total weight of $\dr D$ for the transverse measure
$\alpha$ is equal to the sum of the values of $w$ on the edges of $\gamma$.
If $D\subset N$ (or if $D\subset N'$) then hypothesis (3) of Theorem
\ref{tm:otal} is a direct translation of hypothesis (1) of Theorem
\ref{tm:poly}. We now suppose that it is not the case.

Consider the connected
components $D_1,\cdots, D_N$ of either $(D\cap N)$ or $(D\cap N')$. By the
simplifications done above, each of the $D_i$ has non-empty intersection with
$\dr D$. Since $D$
is a disk, at least two of those connected components, say $D_1$ and $D_2$,
have a boundary which is the union of one (connected) curve of $\dr D$ and of
one connected component of $D\cap N\cap N'$. Then the curve $\dr
D_1\cap \dr D$ (resp. $\dr D_2\cap \dr D$) is in $\dr N\cap \dr \cC$ (or in
$\dr N'\cap \dr \cC$) and begins and ends on a truncation face $f_1$
(resp. $f_2$). Moreover, the deformation made above shows that $\dr D_1\cap
\dr D$ (resp. $\dr D_2\cap \dr D$) does not remain in the faces adjacent to
$f_1$ (resp. $f_2$). Therefore 
hypothesis (2) of Theorem \ref{tm:poly} shows that its total weight for
$\alpha$ is strictly larger than $\pi$. It follows that the total weight of
$\dr D$ for $\alpha$ is strictly larger than $2\pi$, so that hypothesis (3) of
Theorem \ref{tm:otal} is valid.

\medskip
\noindent {\bf 5. Proof that hypothesis (2) of Theorem \ref{tm:otal} is
  satisfied.}
Let $A$ be an essential annulus in $\cC$, such that the total weight of $\dr A$
for $\alpha$ is zero; then each connected component of $\dr A$ is contained in
one of the faces of $\cC$. It follows that each of the connected components of
$\dr A$ enters both $N$ and $N'$, otherwise it would be a closed curve in a
convex cell, and $A$ could not be essential. 

We can suppose that $A$ is transverse to $N\cap N'$. Note that no connected
component of $A\cap (N\cap N')$ is a closed  
non-trivial loop in the interior of $A$, because otherwise --- since $N\cap
N'$ is the disjoint union of a finite number of disks --- $A$ could not be an
essential annulus. 

Moreover, if some connected component of  $A\cap (N\cap
N')$ is a homotopically trivial loops in $A$, we can suppress it by
deforming $A$, as described for disks in the previous step. 
So we can suppose that all the connected components of $A\cap
(N\cap N')$ intersect $\dr A$. The same deformation 
argument can be applied if one
connected component of $A\cap (N\cap N')$ has both endpoints on the same
connected component of $\dr A$. So we can suppose that all connected
components of $A\cap (N\cap N')$ go from one connected component of $\dr A$
to the other. 

Let $D$ be one of the connected components of either $A\cap N$ or $A\cap
N'$. Then $D$ is a strip in $A$, bounded by two curves, each going from one
of the connected components of $\dr A$ to the other. So $D$ is a disk in $N$
(or in $N'$), with boundary in $\dr N$ (or $\dr N'$) a curve which
intersects the closure of at most 
two cells of $\Gamma$ in $\dr N$ (or $\dr N'$)
corresponding to the faces of $\cC$ containing the two boundary curves of
$A$. If $D$ is essential, this is impossible since it contradicts the
hypothesis that the cellular 
decomposition of $\dr' N$ is proper. However, if none of the connected
component of $A\cap N$ or $A\cap N'$ is essential, $A$ is not 
essential. The same argument can be used if $A$ is a Möbius strip, and 
this shows that hypothesis (2) in Theorem
\ref{tm:otal} is valid. 

\medskip
\noindent {\bf 6. End of the proof.}
Theorem \ref{tm:otal} shows that there exists a unique 
hyperbolic metric $h$ on $\cC$
such that $(\cC,h)$ is the convex core of a geometrically finite manifold, with
rational pleating lamination, with pleating lamination given by $\Gamma$ and
transverse measure given by $\alpha$. By the uniqueness, $(\cC,h)$ can be cut
along totally geodesic surfaces to obtain hyperbolic metrics on $N$ and $N'$,
corresponding to truncated hyperideal metrics. Theorem \ref{tm:poly} follows.

\section*{Acknowledgements}

The results presented here owe a lot to a conversation with Francis Bonahon. 
I'm also grateful to Boris Springborn for some conversations
directly related to this paper.

\bibliographystyle{alpha}

\end{document}